\documentclass{article}

\usepackage{amssymb}
\usepackage{amsmath}
\usepackage{color}
\usepackage{graphicx}
\usepackage{subcaption}
\usepackage{geometry} 
\usepackage[english,francais]{babel}
\def\R{\mathbb{R}}
\def\Z{\mathbb{Z}}

\def\bu{\mathbf{u}}
\def\bv{\mathbf{v}}
\def\be{\mathbf{e}}
\def\cA{\mathcal{A}}
\def\cG{\mathcal{G}}
\def\dsp{\displaystyle}

\newtheorem{theorem}{Theorem}[section]
\newtheorem{lemma}[theorem]{Lemma}

\newtheorem{e-proposition}[theorem]{Proposition}

\newtheorem{e-definition}[theorem]{Definition\rm}

\def\og{\leavevmode\raise.3ex\hbox{$\scriptscriptstyle\langle\!\langle$~}}
\def\fg{\leavevmode\raise.3ex\hbox{~$\!\scriptscriptstyle\,\rangle\!\rangle$}}

\begin{document}
\centerline{}


\selectlanguage{english}

\begin{center}
{\Large{Existence of guided waves due to a lineic perturbation of a 3D periodic medium }}\\[0.5cm]

{\large{B\'erang\`ere Delourme$^{a,b}$, Patrick Joly$^{b}$, Elizaveta Vasilevskaya$^{a,b}$}}\\[0.5cm]

{\small $a$: Universit\'e Paris 13, Sorbonne Paris Cit\'e, LAGA, UMR
  7539, 93430 Villetaneuse, France} \\
{\small $b$: POEMS, UMR 7231 CNRS/ENSTA/Inria, ENSTA ParisTech, 828 boulevard des
Mar\'echaux, 91762 Palaiseau Cedex, France}\\
\end{center}






\begin{abstract}
\selectlanguage{english}
In this note, we exhibit a three dimensional structure that permits to guide waves.  This structure is obtained by a geometrical perturbation of a 3D periodic domain that consists of a three dimensional grating of equi-spaced thin pipes oriented along three orthogonal directions. Homogeneous Neumann boundary conditions are imposed on the boundary of the domain.  The diameter of the section of the pipes, of order $\varepsilon>0$, is supposed to be small. We prove that, for $\varepsilon$ small enough, shrinking the section of one line of the grating by a factor of $\sqrt{\mu}$ ($0<\mu<1$) creates guided modes that propagate along the perturbed line. Our result relies on the asymptotic analysis (with respect to $\varepsilon$) of the spectrum of the Laplace-Neumann operator in this structure. Indeed, as $\varepsilon$ tends to $0$, the domain tends to a periodic graph,  and the spectrum of the associated limit operator can be computed explicitly. \\[2ex]
\end{abstract}



\selectlanguage{english}

\section{Statement of the problem}\label{}

\noindent Let $\omega_1$, $\omega_2$ and $\omega_3$ be three Lipschitz bounded domains of $\R^2$ of same area ($|\omega_1| = |\omega_2| =|\omega_3|$)  containing the origin $(0,0)$,  let $\varepsilon>0$ be a parameter (that is going to be small), 
and let  $a_1$, $a_2$ and $a_3$ be three positive real numbers. We denote by $(\be_i)_{i \in \{1, 2 , 3\}}$, the standard basis of $\R^3$. For any $(k, \ell) \in \Z^2$, we consider the three dimensional domain $D_{k,\ell,3}^\varepsilon$ defined by  
$$
D_{k,\ell,3}^\varepsilon = \left\{ (x_1 , x_2 , x_3) \in \R^3 \mbox{ such that }  \big(({x_1 - a_1 k})/{\varepsilon},(x_2 -a_2 \ell)/{\varepsilon}\big)  \in \omega_3 \right\},  
$$
which is an unbounded cylinder of constant cross section $\varepsilon \omega_3$. It is infinite along the $\be_3$ direction (invariant with respect to $x_3$) and contains the  point $(a_1 k, a_2 \ell ,0)$.  
Similarly, for any $(k, \ell) \in \Z^2$, we define the domains 
\begin{align*}
& D_{k,\ell,1}^\varepsilon =   \left\{ (x_1 , x_2 , x_3) \in \R^3 \mbox{ such that }  \big(({x_2 - a_2 k})/{\varepsilon},(x_3 -a_3 \ell)/{\varepsilon}\big)  \in \omega_1 \right\},\\
& D_{k,\ell,2}^\varepsilon =   \left\{ (x_1 , x_2 , x_3) \in \R^3 \mbox{ such that }  \big(({x_1 - a_1 k})/{\varepsilon},(x_3 -a_3 \ell)/{\varepsilon}\big)  \in \omega_2 \right\},
\end{align*}
and we consider the periodic domain $\Omega_\varepsilon$  given by 
\begin{equation}
\Omega_\varepsilon = \bigcup_{i\in \{1,2,3\}} \bigcup_{(k,\ell)\in \Z^2}  D_{k,\ell,i}^\varepsilon. 
\end{equation}
\noindent The domain $\Omega_\varepsilon$ is a three dimensional grating of equi-spaced parallel pipes (of constant cross section) oriented along the three orthogonal directions $\be_1$, $\be_2$ and $\be_3$. It is $a_j$-periodic with respect to $x_j, j=1,2,3$. Moreover, the points  $(k a_1,  \ell a_2, m a_3)$, $(k,\ell,m) \in \Z^3$, belong to $\Omega_\varepsilon$. \\
\begin{figure}[htbp]
        \centering
        \begin{subfigure}[b]{0.4\textwidth}
      \begin{center}
                \includegraphics[width=\textwidth,trim=0cm 0cm 0cm 0.5cm, clip]{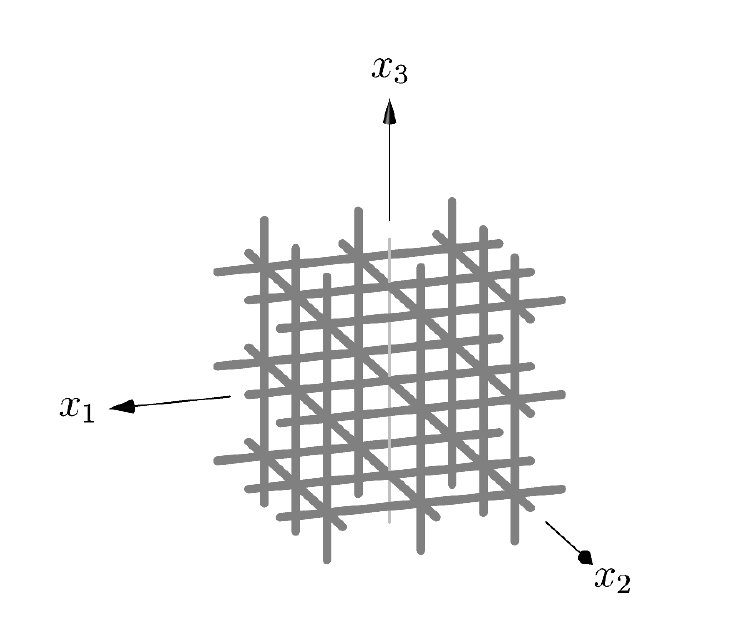}
                \caption{The perturbed domain $\Omega_\varepsilon^\mu$}
                \label{fig:OmegaPerturbe}
                \end{center}
        \end{subfigure}
       \quad
        \begin{subfigure}[b]{0.4\textwidth}
        \begin{center}
                \includegraphics[width=\textwidth, trim=0cm -1cm 0cm 0.5cm, clip]{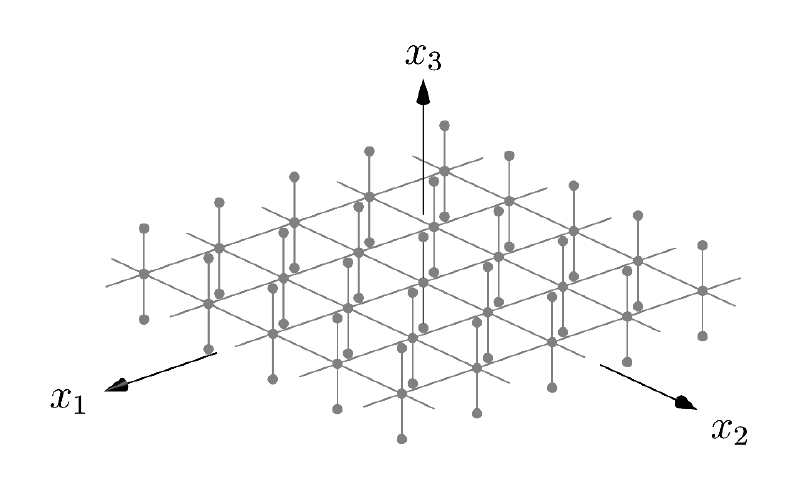}
                \caption{The graph $\mathcal{G}$.}
                \label{fig:Graphe}
                \end{center}
        \end{subfigure}
         \caption{Illustration of the perturbed periodic domain $\Omega_\varepsilon^\mu$ and the limit graph  $\mathcal{G}$.}\label{fig:dessinsDomaines}
\end{figure}

\noindent
\noindent In order to create guided modes, we introduce a linear defect (see~\cite{AmmariSantosa}-\cite{KuchmentOng}-\cite{Sonia}) in the periodic structure by modifying the section size of one pipe of the grating (it is conjectured that guided modes cannot appear in the purely periodic structure, see~\cite{Friedlander} for the proof in the case of a symmetric medium).
 More precisely, we assume that the domain $D_{0,0,3}^\varepsilon$ is replaced with the domain  
$$
D_{0,0,3}^{\varepsilon,\mu} =  \left\{ (x_1 , x_2 , x_3) \in \R^3 \mbox{ such that }  \big(x_1 /{(\sqrt{\mu}\varepsilon)},x_2 /{(\sqrt{\mu}\varepsilon)}\big)  \in \omega_3 \right\},  
$$
where $\mu$ is a positive parameter. In other words, we enlarge ($\mu>1$) or shrink ($0< \mu <1$) the section of one pipe of the domain by a factor $\mu$ (see Fig~\ref{fig:OmegaPerturbe}). The corresponding  perturbed domain is denoted by  $\Omega_\varepsilon^\mu$. Its precise definition is given by
\begin{equation}
\Omega_\varepsilon^\mu = \left( \bigcup_{i\in \{1,2\}} \bigcup_{(k,\ell)\in \Z^2}  D_{k,\ell,i}^\varepsilon \right) \bigcup  \left( \bigcup_{(k,\ell)\in \Z^2\setminus \{(0,0) \}} D_{k,\ell,3}^\varepsilon \right) \bigcup D_{0,0,3}^{\varepsilon,\mu}.
\end{equation}
$\Omega_\varepsilon^\mu$ is still $a_3$-periodic with respect to $x_3$. However,  the presence of the perturbed pipe $D_{0,0,0}^{\varepsilon,\mu}$ breaks the periodicity with respect to $x_1$ and $x_2$. We emphasize that the domain $\Omega_\varepsilon^\mu$ (as well as $\Omega_\varepsilon$)  tends to a 3D periodic graph as $\varepsilon$ tends to $0$. \\

\noindent We look for guided modes, i.e. solutions to the wave equation $\partial_t^2u - \Delta u = 0$ in $\Omega_\varepsilon^\mu$, satisfying homogeneous Neumann boundary conditions on $\partial \Omega_\varepsilon^\mu$ (see~\cite{NazarovDirichletLadder} for the  investigation of the Dirichlet case), that propagate along the defect pipe $D_{0,0,0}^{\varepsilon,\mu}$ (i.e. in the $\be_3$ direction) but stay confined in the transversal directions. More precisely, denoting by $B_\varepsilon^\mu$ the restriction of the domain $\Omega_\varepsilon^\mu$ to the band $|x_3|<a_3/2$, 
\begin{equation}\label{Bepsilonmu}
B_\varepsilon^\mu = \left\{ (x_1, x_2, x_3) \in \Omega_\varepsilon^\mu \; \mbox{such that} \; |x_3| < {a_3}/{2}\right\},
\end{equation} 
we search solutions of the form $
u(x_1,x_2,x_3,t)=v(x_1,x_2,x_3)e^{i \omega t -\beta x_3}$, 
where $\beta$ is a real parameter and $v(x_1,x_2,x_3)\in L_2(B_\varepsilon^\mu)$ is an $a_3$-periodic function in $x_3$.  In fact, it is easily seen that the $\beta$-quasiperiodic fonction $v(x_1,x_2,x_3)e^{-i\beta x_3}$  is an eigenfunction of the operator 
\begin{equation}\label{DefinitionAvarepsilonmu}
A_{\varepsilon}^{\mu}(\beta): D\left(A_{\varepsilon}^{\mu}(\beta)\right) \subset L_2\left(B_{\varepsilon}^{\mu}\right)\rightarrow L_2\left(B_{\varepsilon}^{\mu}\right), \quad A_{\varepsilon}^{\mu}(\beta)=-\Delta u \;\; \mbox{in} \; B_{\varepsilon}^{\mu},\\
\end{equation}
with
$
D\left(A_{\varepsilon}^{\mu}(\beta)\right)=\big\{u\in H^1_{\Delta}\left(B_{\varepsilon}^{\mu}\right),
\displaystyle
\left.u\right|_{\Sigma^+}=e^{-i\beta}\left.u\right|_{\Sigma^-}, \displaystyle\left.\partial_{x_3} u \right|_{\Sigma^+}=e^{-i\beta}\displaystyle\left. \partial_{x_3} u \right|_{\Sigma^-},\; \left.\partial_n u \right|_
{\partial B_{\varepsilon}^{\mu}\setminus\Sigma^\pm}\!\!=0\big\},
$
where 
$$H^1_{\Delta}(B_{\varepsilon}^{\mu}) = \left\{ u \in H^1(B_{\varepsilon}^{\mu}), \mbox{s.t.} \; \Delta u \in L_2(B_{\varepsilon}^{\mu}) \right\} \mbox{  and } \Sigma^\pm = \left\{ (x_1, x_2, x_3) \in \partial B_{\varepsilon}^{\mu}, x_3 = \pm a_3/2\right\}.$$
\noindent To study the spectral properties of $A_{\varepsilon}^{\mu}(\beta)$, we investigate its (formal) limit $\cA^\mu(\beta)$ as $\varepsilon$ tends to $0$.  The operator $\cA^\mu(\beta)$ is defined on the limit graph $\cG$ (see Fig.~\ref{fig:Graphe}) and its spectrum can be explicitly computed. In particular, its spectrum has infinitely many gaps (Lemma~\ref{LemmeSpectreEssentiel}), i.e. open intervals $(a,b) \subset \R$ such that the intersection of $[a,b]$ with the spectrum is reduced to $\{a,b\}$.
 Moreover, for $\mu<1$, there is at least one eigenvalue in each gap (Lemma~\ref{LemmeValeursPropresGraphe}). Since, in addition, for $\varepsilon>0$ sufficiently small,  the spectrum of $A_{\varepsilon}^{\mu}(\beta)$ is close to the spectrum of $\cA^{\mu}(\beta)$, the existence of guided modes is guaranteed (Theorem~\ref{TheoremeConvergence}).  
\section{The spectrum  of the limit operator $\cA^\mu(\beta)$}
\subsection{Definition of the limit operator $\cA^\mu(\beta)$}
\noindent The limit operator $\cA^\mu(\beta)$ is defined on the infinite periodic graph $\cG=\bigcap_{\varepsilon>0} B_{\varepsilon}^{\mu}$  obtained as the limit of $B_{\varepsilon}^{\mu}$ as $\varepsilon$ tends to $0$:
 $\cG$ is made of the vertices $\{ v_{k,\ell} = (ka_1, \ell a_2, 0), v_{k,\ell}^\pm = (ka_1 , \ell a_2, \pm a_3/2), (k, \ell) \in \Z^2 \}$ connected by the edges 
 $$\{ e_{k+1/2, \ell} = (v_{k,\ell}, v_{k+1,\ell}), \;\;  e_{k, \ell+1/2} = (v_{k,\ell}, v_{k,\ell+1}),  \;\; e_{k, \ell}^\pm =(v_{k,\ell}, v_{k,\ell}^\pm), \;\; (k, \ell) \in \Z^2  \}.$$
It is $a_1$-periodic with respect to $x_1$ and  $a_2$-periodic with respect to $x_2$ (see Fig.~\ref{fig:Graphe}). \\

\noindent For any function $u$  defined on $\cG$, we denote by $\bu_{k, \ell}$ (resp. $\bu_{k, \ell}^\pm$) its value at the vertex $v_{k,\ell}$ (resp. $v_{k,\ell}^\pm$). The restriction of $u$ to the edge $e_{k+1/2,\ell}$ (resp. $e_{k,\ell+1/2}$ and $e_{k,\ell}^\pm$) is denoted by $u_{k+1/2,\ell}(x_1)$ (resp.  $u_{k,\ell+1/2}(x_2)$ and $u^\pm_{k,\ell}(x_3)$).\\

\noindent The definition of $\cA^\mu(\beta)$  also requires the introduction of the function spaces $L_2^\mu(\cG)$ and  $H^2(\cG)$ defined as 
\begin{equation}\label{DefL2muH2mu}
L_{2}^{\mu}\left(\cG\right)= \left\{u:\|u\|_{L_{2}^{\mu}\left(\cG\right)} < + \infty\right\}, \quad H^{2}\left(\cG\right)= \left\{u \in C(\cG):\|u\|_{H^{2}\left(\cG\right)} <  + \infty \right\}, 
\end{equation}
where,
\begin{equation}\label{DefNormeL2mu}
\|u\|^2_{L_{2}^{\mu}\left(\cG\right)}  = \sum\limits_{(k,\ell)\in \mathbb{Z}^2} \Big({w^{\mu}_{k\ell}}
\sum_\pm \|u_{k,\ell}^\pm\|_{L_2(e_{k,\ell}^\pm)}^2
+\|u_{k+\frac{1}{2},\ell}\|_{L_2(e_{k+\frac{1}{2},\ell})}^2+\|u_{k,\ell+\frac{1}{2}}\|_{L_2(e_{k,\ell +\frac{1}{2}})}^2 \; \Big),
\end{equation}
\begin{equation}\label{DefNormeH2mu}
\|u\|^2_{H^{2}\left(\cG\right)}  = \sum\limits_{(k,\ell)\in \mathbb{Z}^2} \Big( 
\sum_\pm \|u_{k,\ell}^\pm\|_{H^2(e_{k,\ell}^\pm)}^2
+\|u_{k+\frac{1}{2},\ell}\|_{H^2(e_{k+\frac{1}{2},\ell})}^2+\|u_{k,\ell+\frac{1}{2}}\|_{H^2(e_{k,\ell +\frac{1}{2}})}^2 \; \Big),
\end{equation}
and, for any $(k, \ell) \in \Z^2$, $w^{\mu}_{k,\ell}$ is the weight coefficient  equal to $\mu$ for $k = \ell =0$ and $1$ otherwise.\\

\noindent The unbounded limit operator in $L_{2}^{\mu}(\cG)$  has domain
\begin{multline}
\label{Kirch_3d}
D\left(\cA^{\mu}(\beta)\right)=\left\{u\in H^2\left(\cG\right):\quad\forall (k,\ell) \in\mathbb{Z}^2,\quad
\vphantom{\left. u'_{j+\frac{1}{2}}(0)-u'_{j-\frac{1}{2}}(1)+\mu(u_j^{+})'(0)-\mu(u_j^{-})'(0)=0,\qquad\: \forall j\in\mathbb{Z}\right\}}\right.
 \bu_{k,\ell}^+=e^{-i\beta} \bu_{k,\ell}^-,\quad (u^{+}_{k,\ell})'\left(\textstyle{a_3}/{2}\right)= e^{-i\beta}(u^{-}_{k,\ell})'\left(-\textstyle{a_3}/{2}\right),\\
\qquad\;\left. u'_{k+\frac{1}{2},\ell}(k a_1)-u'_{k-\frac{1}{2},\ell}(k a_1)+u'_{k,\ell+\frac{1}{2}}(\ell a_2)-u'_{k,\ell-\frac{1}{2}}(\ell a_2)+{w^{\mu}_{k,\ell}}\left((u_{k,\ell}^{+})'(0)-(u_{k,\ell}^{-})'(0)\right)=0\right\},
\end{multline}
and  is defined by
\begin{equation}
\forall \; u \in D\left(\cA^{\mu}(\beta)\right), \quad \cA^\mu(\beta) u = - u'' \quad  \mbox{on any edge of the graph} \, \cG.  
\end{equation}
The functions of $D\left(\cA^{\mu}(\beta)\right)$ are continuous on $\cG$ and $\beta$ quasi-periodic. Moreover, they satisfy the Kirchhoff conditions~\eqref{Kirch_3d} that enforce the weighted sum of the outward derivatives of $u$ to vanish at each vertex $v_{k,\ell}$ ($(k, \ell) \in \Z^2$).  We point out that the perturbation, which results from a geometrical modification of the domain for the problem~\eqref{DefinitionAvarepsilonmu}, is taken into account at the limit by means of the Kirchhoff condition~\eqref{Kirch_3d} at the vertex $v_{0,0}$ ($w^{\mu}_{0,0}=\mu$). The formal derivation of the limit model can be found in~\cite{KuchmentZeng}. 
\noindent It is easily verified that the operator $\cA^{\mu}(\beta)$ is self-adjoint (for the weighted scalar product associated with~\eqref{DefNormeL2mu}), see also~\cite{KuchmentSurvey}. The objective of the following two sections is to study the spectrum of $\cA^{\mu}(\beta)$.
\subsection{Characterization and properties of the essential spectrum of $\cA^{\mu}(\beta)$}
\noindent By a compact perturbation argument, one can prove that $\sigma_{ess}(\cA^{\mu}(\beta))= \sigma(\cA(\beta))$, where $\cA(\beta)=\cA^{1}(\beta)$ is the purely periodic operator corresponding to  $\cA^{\mu}(\beta)$ for $\mu=1$. The computation of its spectrum relies on the Floquet-Bloch theory (see~\cite{ReedSimon}). More precisely, we can prove that 
 $\lambda = \omega^2 \in \sigma(\cA(\beta))$ if and only if either $\omega =0$ and $\beta=0$ or $\omega \neq 0$ and there exists $(k_1,k_2) \in [0,\pi]^2$ such that
\begin{multline}
\label{condition_spectre_3d}
\sin{(\omega a_2)}\sin{(\omega a_3)}\left(\cos{(\omega a_1)-\cos{k_1}}\right)+
\sin{(\omega a_3)}\sin{(\omega a_1)}\left(\cos{(\omega a_2)-\cos{k_2}}\right)\\
+\sin{(\omega a_1)}\sin{(\omega a_2)}\left(\cos{(\omega a_3)-\cos{\beta}}\right)=0.
\end{multline} 
Based on the previous characterization, we prove that the operator $\cA(\beta)$ has a countable infinity of gaps that can be separated into three categories (see~\cite{TheseElizaveta} for the proof): \\
\begin{lemma}\label{LemmeSpectreEssentiel} The following properties hold :
\begin{enumerate}
\item[1-] $\sigma_{1} \cup \sigma_2 \cup \sigma_3 \subset \sigma(\cA(\beta))$, where 
$$
\sigma_{i}  = \left\{ \left({\pi n}/{a_i}\right)^2, n \in \Z  \right\} \; \mbox{for} \; i \in \{ 1,2\},  \; \mbox{and} \; \sigma_{3}  = \left\{ \left((\pm \beta + 2 \pi n)/a_3\right)^2, n \in \Z  \right\}.
$$
\item[2-] For any $\beta \in [0, \pi]$, the operator $\cA(\beta)$ has infinitely many gaps whose ends tend to infinity.
\item[3-] Let 
$
\mathcal{W}(\beta) =
\left\{{\pi n}/{a_3},\; n\in\mathbb{N}^*\right\} \mbox{if } \beta\notin\{0,\pi\}$ and
$ \mathcal{W}(\beta) = \left\{{\beta}/{a_3}+{(2n+1)\pi}/{a_3},\; n\in\mathbb{N}^*\right\} \mbox{if }  \beta\in\{0,\pi\},
$.\\[4pt]
If an interval $(\omega_b^2, \omega_t^2)$ is a spectral gap of $\cA(\beta)$, then, one of the following possibilities holds: \\[-2pt]
\begin{enumerate}
\item[(i)] $\omega_b^2 \notin \sigma_1 \cup \sigma_2$,  $\omega_t^2 \notin \sigma_1 \cup \sigma_2$, and there is a unique $\omega_0 \in (\omega_b, \omega_t) \cap \mathcal{W}(\beta)$. 
\item[(ii)]  $\omega_b^2 \in \sigma_1 \cup \sigma_2$, $\omega_t^2 \notin \sigma_1 \cup \sigma_2$ and $\mathcal{W}(\beta)\cap (\omega_b,\omega_t) = \varnothing$. 
\item[(iii)]  $\omega_b^2 \notin \sigma_1 \cup \sigma_2$, $\omega_t^2 \in \sigma_1 \cup \sigma_2$ and $\mathcal{W}(\beta)\cap (\omega_b,\omega_t) = \varnothing$. 
\end{enumerate}
\end{enumerate}
\end{lemma} 

\subsection{Computation of the discrete spectrum}

\noindent Let us now determine the discrete spectrum of $\cA^\mu(\beta)$.  If $\lambda=\omega^2$ is an eigenvalue of $\cA^\mu(\beta)$, then the corresponding eigenfunction $u \in D(\cA^\mu(\beta))$ satisfies the linear differential equation $u'' + \omega^2 u=0$ on each edge of the graph $\mathcal{G}$.  Solving explicitly this equation (on each edge), taking into account the quasi-periodicity of $u$ and the Kirchhoff conditions~\eqref{Kirch_3d}, we can replace the eigenvalue problem $\cA^\mu(\beta) u = \lambda u$ with a set of finite differences equations for $(\bu_{\ell, k})_{(\ell,k)\in \Z^2}$:\\[-2ex]

\begin{lemma}\label{LemmeDifferencesFinies} Assume that $\omega \notin \left\{\pi\Z/a_1\right\}\cup\left\{\pi\Z/a_2\right\}\cup\left\{\pi \Z/a_3\right\}$. $u \in D(\cA^\mu(\beta))$ is an eigenfunction of $\cA^\mu(\beta)$ if and only if $(\bu_{k,\ell})_{(k,\ell)\in \Z^2}$ belongs to $\ell_2(\Z^2)$ and satisfies
\begin{equation}\label{FDkl}
\left\{ \begin{array}{l}\dsp \frac{\bu_{k+1,\ell}}{\sin{(\omega a_1)}}+\frac{\bu_{k-1,\ell}}{\sin{(\omega a_1)}}+
\frac{\bu_{k,\ell+1}}{\sin{(\omega a_2)}}+\frac{\bu_{k,\ell-1}}{\sin{(\omega a_2)}}-
2g_{\beta}(\omega)\bu_{k,l}=0,\qquad \forall  (k,\ell)\in\mathbb{Z}^2 \setminus \left\{ (0,0) \right\},\\[12pt]
\dsp \frac{\bu_{1,0}}{\sin{(\omega a_1)}}+\frac{\bu_{-1,0}}{\sin{(\omega a_1)}}+
\frac{\bu_{0,1}}{\sin{(\omega a_2)}}+\frac{\bu_{0,-1}}{\sin{(\omega a_2)}}-
2g_{\beta}(\omega)\bu_{0,0}=2(\mu-1)\frac{\cos{(\omega a_3)-\cos{\beta}}}{\sin{(\omega a_3)}}\bu_{0,0}, 
\end{array} \right.\end{equation}
where we have defined
$
\dsp g_{\beta}(\omega)=\frac{1}{\tan{(\omega a_1)}}+\frac{1}{\tan{(\omega a_2)}}+\frac{\cos{(\omega a_3)}-\cos{\beta}}{\sin{(\omega a_3)}} \; .
$
\end{lemma} 
~\\
As well-known, finite difference schemes may be investigated using the discrete Fourier transform
\begin{equation}\label{DFT}
\mathcal{F} : 
\bv = (\bv_{k, \ell})_{(k, \ell) \in \Z^2} \mapsto \mathcal{F}(\bv) =   \widehat{\bv}, \quad \widehat{\bv}(\xi,\eta)=\sum\limits_{k,\ell\in\mathbb{Z}}e^{i(k\xi+\ell\eta)}\bv_{k,\ell},\qquad (\xi,\eta)\in[0,2\pi]^2,
\end{equation}
where $\mathcal{F} $ is an isometry between $\ell_2(\Z^2)$ and $L_2((0, 2\pi)^2)$. This, together with Lemma~\ref{LemmeDifferencesFinies}, provides the following characterization for the discrete spectrum of $\cA^\mu(\beta)$:\\[-2ex]
\begin{lemma}\label{LemmaDFT}Assume that $\omega \notin \left\{\pi\Z/a_1\right\}\cup\left\{\pi\Z/a_2\right\}\cup\left\{\pi \Z/a_3\right\}$. $u \in D(\cA^\mu(\beta))$ is an eigenfunction of $\cA^\mu(\beta)$ if and only if the discrete Fourier transform $\widehat{\bu}$ of $(\bu_{\ell, k})_{(\ell,k)\in \Z^2}$ belongs to $L_2((0,2\pi)^2)$ and satisfies
\begin{equation}
\label{Equchapeau}
\left(  f(\xi, \eta, \omega)-\phi_{\beta}(\omega)\right) \;  \widehat{\bu}(\xi,\eta) = (\mu-1) \phi_{\beta}(\omega)  {\bu}_{0,0}.
\end{equation}
where 
$
\dsp \phi_{\beta}(\omega)  =  \frac{\cos{(\omega a_3)-\cos{\beta}}}{\sin{(\omega a_3)}}  \; \mbox{and} \; f(\xi, \eta, \omega) =  \frac{\cos{\xi}-\cos{(\omega a_1)}}{\sin{(\omega a_1)}}+
\frac{\cos{\eta}-\cos{(\omega a_2)}}{\sin{(\omega a_2)}}.
$
\end{lemma} 
~\\[2pt]
Under the assumption of Lemma~\ref{LemmaDFT}, \eqref{condition_spectre_3d} can be written as $f(\xi, \eta, \omega)-\phi_{\beta}(\omega)=0$. It follows that  $\lambda = \omega^2$ does not belong to $\sigma_{ess}(\cA^\mu(\beta))$ if and only if,  for any $(\xi, \eta) \in [0,2 \pi]^2$, $\phi_{\beta}(\omega) - f(\xi, \eta, \omega)$ does not vanish. As a consequence,  as soon as $\lambda  = \omega^2 \notin \sigma_{ess}(\cA^\mu(\beta))$, the function $(\xi, \eta)\mapsto \phi_{\beta}(\omega)/( \phi_{\beta}(\omega) - f(\xi, \eta, \omega))$ is continuous and bounded. Then, the inverse discrete Fourier transform can be applied to~\eqref{Equchapeau} to obtain
$$
\bu_{k,\ell}=\frac{(1-\mu)\bu_{0,0}}{4\pi^2}\int\limits_{(0,2\pi)^2} \frac{\phi_{\beta}(\omega)}{\phi_{\beta}(\omega)-f(\xi,\eta,\omega)}e^{-i(k\xi+\ell\eta)}d\xi d\eta,\quad  \forall \, (k,\ell) \in\mathbb{Z}^2.
$$
Writing the previous relation for $k = \ell=0$ yields the following criterion of existence of an  eigenvalue:\\[-2ex]
\begin{lemma}
Assume that $\omega \notin \left\{\pi\Z/a_1\right\}\cup\left\{\pi\Z/a_2\right\}\cup\left\{\pi \Z/a_3\right\}$ and that $\lambda = \omega^2 \notin \sigma_{ess}(\cA^\mu(\beta))$. Then, $\lambda$ is an eigenvalue of $\cA^\mu(\beta)$ if and only if
\begin{equation}\label{equationVPFinale}
\mu = 1 - F_\beta(\omega) \quad \mbox{where} \quad F_\beta(\omega)  =\Big( \; \frac{1}{4 \pi^2} \int\limits_{(0,2 \pi)^2}\frac{\phi_\beta(\omega)} {\phi_{\beta}(\omega)-f(\xi,\eta, \omega)} d\xi d\eta \; \Big)^{-1}.
\end{equation} 
\end{lemma}   
The study of the behavior of the function $F_\beta$ leads to the existence of at least one eigenvalue in each gap of $\cA^\mu(\beta)$ as soon as $\mu<1$, the minimal number of eigenvalues in each gap depending on the type of gaps (cf. Lemma.~\ref{LemmeSpectreEssentiel}-3 for the classification):\\[2pt]
\begin{lemma}\label{LemmeValeursPropresGraphe}
For $\mu>1$, the operator $\cA^\mu(\beta)$ has no eigenvalue. For $0<\mu<1$, let $(\omega_b^2, \omega_t^2)$ be a spectral gap of the operator $\cA^\mu(\beta)$: 
\begin{enumerate}
\item[(a)] If $(\omega_b^2, \omega_t^2)$ is a gap of type (i), then $\cA^\mu(\beta)$  has at least two eigenvalues $\lambda_1 = \omega_1^2$ and $\lambda_2= \omega_2^2$ that satisfy
$\omega_b < \omega_1 < \omega_0 < \omega_2< \omega_t$ (see Lemma.~\ref{LemmeSpectreEssentiel}-3 for the definition of $\omega_0$).
\item[(b)] If $(\omega_b^2, \omega_t^2)$ is a gap of type (ii) or (iii), then $\cA^\mu(\beta)$  has at least one eigenvalue $\lambda_1 = \omega_1^2$ such that $\omega_b < \omega_1 < \omega_t$.  
\end{enumerate}   
\end{lemma}
~\\[2pt]
The sketch of the proof of the previous lemma is the following, a complete proof being available in \cite{TheseElizaveta} (Theorem 5.2.1): First, one can verify that $F_\beta(\omega)\geq0$ in any gap, which, together with~\eqref{equationVPFinale} proves that $\cA^\mu(\beta)$ has no eigenvalue for $\mu>1$. 
 Then, if $(\omega_b^2, \omega_t^2)$ is a gap of type (i), one can show that 
 $$
 \lim_{\omega\rightarrow \omega_b^+} (1-F_\beta(\omega)) =  \lim_{\omega\rightarrow \omega_t^-} (1-F_\beta(\omega)) =1
 \mbox{  
 and that }\lim_{\omega\rightarrow \omega_0} (1-F_\beta(\omega)) = 0.$$  By continuity of $F_\beta$ inside the gap,  (a) directly results from the intermediate value theorem and~\eqref{equationVPFinale}. 
If $(\omega_b^2, \omega_t^2)$ is a gap of type (ii), the intermediate value theorem also permits us to conclude since  $$\lim_{\omega \rightarrow \omega_b} (1 -F_\beta(\omega)) \leq 0 \mbox{ and } \lim_{\omega\rightarrow \omega_t^-} (1-F_\beta(\omega)) =1. $$ 
A similar argument works for a gap of type (iii).

\section{Guided modes for the operator $A_\varepsilon^\mu(\beta)$: an asymptotic result}
\noindent Finally, thanks to the general result~\cite{Post} (Theorem 2.13 convergence of the spectrum of $A_\varepsilon^\mu(\beta)$ toward the spectrum of ${\cA}^{\mu}(\beta)$), we can prove the following result of existence of eigenvalue for the operator $A_\varepsilon^\mu(\beta)$:\\[2pt]
\begin{theorem}
\label{TheoremeConvergence}
Let $\mu \in (0,1)$,  $(\lambda_b, \lambda_t)$ be a spectral gap of the operator ${\cA}^{\mu}(\beta)$  and $\lambda_0 \in (\lambda_b,\lambda_t)$ be a (simple) eigenvalue of this operator. Then, there exists $\varepsilon_0 >0$ such that if $\varepsilon<\varepsilon_0$ the operator $A_\varepsilon^{\mu}(\beta)$ has an eigenvalue $\lambda_{\varepsilon}$ inside a spectral gap $(\lambda_b^\varepsilon,\lambda_t^\varepsilon )$. Moreover, 
$
\lambda_{\varepsilon}=\lambda_0+O(\sqrt{\varepsilon}). 
$
\end{theorem} 



\end{document}